\documentclass[12pt, requno]{amsart}
\usepackage{amsmath}
\usepackage{amssymb}
\usepackage{epsfig}
\usepackage{graphicx}
\usepackage{color}
\definecolor{shadecolor}{gray}{0.875}
\usepackage{amscd}
\usepackage{comment}
\usepackage{hyperref}
\usepackage{tikz-cd}

\numberwithin{equation}{section}

\input xy
\xyoption{all}

\calclayout
\allowdisplaybreaks[3]

\theoremstyle{plain}
\newtheorem{prop}{Proposition}[]

\newtheorem{theo}[prop]{Theorem}
\newtheorem{coro}[prop]{Corollary}

\newtheorem{lemm}[prop]{Lemma}

\theoremstyle{definition}
\newtheorem{defi}[prop]{Definition}

\def\lra{\longrightarrow}

\def\cC{{\mathcal C}}

\def\cE{{\mathcal E}}

\def\cL{{\mathcal L}}

\def\cN{{\mathcal N}}

\def\cU{{\mathcal U}}
\def\cV{{\mathcal V}}
\def\cW{{\mathcal W}}
\def\cX{{\mathcal X}}
\def\cY{{\mathcal Y}}

\def\bC{{\mathbb C}}
\def\bD{{\mathbb D}}
\def\bF{{\mathbb F}}

\def\bP{{\mathbb P}}

\def\bZ{{\mathbb Z}}

\def\bk{\Bbbk}

\def\rB{{\mathrm B}}
\def\rH{{\mathrm H}}

\def\rP{{\mathrm P}}

\def\Pic{\mathrm{Pic}}

\def\Hilb{\mathrm{Hilb}}

\def\Spec{\mathrm{Spec}}

\def\Pic{\mathrm{Pic}}

\def\Sect{\mathrm{Sect}}

\makeatother
\makeatletter

\author{Sho Tanimoto}
\address{Graduate School of Mathematics, Nagoya University, Furocho Chikusa-ku, Nagoya, 464-8602, Japan}
\email{sho.tanimoto@math.nagoya-u.ac.jp}

\author[Yu. Tschinkel]{Yuri Tschinkel}
\address{
  Courant Institute,
  251 Mercer Street,
  New York, NY 10012, USA
}
\email{tschinkel@cims.nyu.edu}

\address{Simons Foundation\\
160 Fifth Avenue\\
New York, NY 10010\\
USA}

\title[Homological stability and weak approximation]{Homological stability\\ and weak approximation}

\begin{document}

\maketitle

\begin{abstract}
We investigate homological stability for the space of sections of Fano fibrations over curves in the context of weak approximation, and establish it for projective bundles, as well as for conic and quadric surface bundles over curves.
\end{abstract}

\tableofcontents

\section{Introduction}

Let $B$ be a smooth projective geometrically irreducible curve over a field $\bk$ and  
$F=\bk(B)$ its function field. Let $X$ be a smooth projective variety over $F$  and   
\begin{equation}
    \label{eqn:fibr}
\pi:\cX\to B
\end{equation}
its smooth integral model, i.e., a flat proper morphism from a smooth projective $\cX$ over $\bk$, with generic fiber $X$. 
Our point of departure is the connection between arithmetic properties of $X$ over $F$ and geometric properties of  
spaces of sections of $\pi$, over $\bk$. Of particular interest are cases when $\bk=\bF_q$, a finite field, or $\bk=\bC$. 

Concretely, let $\omega^{-1}_\pi$ be the relative anticanonical class, assumed to be 
ample on the generic fiber $X$. 
We are interested in understanding the space of sections of $\pi$
$$
\Sect(\cX/B,h):=\{ \sigma:B\to \cX \} 
$$
of height 
\begin{equation*}  
\label{eqn:height}
h:=\deg(\sigma^*\omega_\pi^{-1}).  
\end{equation*}
These spaces have been studied from different but related perspectives: 
\begin{itemize}
\item {\em Manin's conjecture}, see, e.g., \cite{Bat88},\cite{FMT}, \cite{BM}, \cite{Peyre}, \cite{BT}, \cite{Peyre03}, \cite{Peyre17}, \cite{LST18},
\item {\em Homological stability} \cite{CJS00}, \cite{DT24}, \cite{DLTT25}. 
\end{itemize}
Manin's conjecture concerns {\em asymptotics} of points (sections) of bounded height. When $\bk=\bF_q$, this translates into understanding the growth of  
\begin{equation*} \label{eqn:sect}
\# \Sect(\cX/B,h)(\bF_q),  \quad h\to\infty. 
\end{equation*}
In turn, these numbers can be accessed via Grothendieck’s Lefschetz trace formula, which motivates the study of {\em topology} of complex points of spaces of sections of fibrations \eqref{eqn:fibr} defined over $\bk=\bC$. 
Homological stability asserts stabilization of homology of these spaces, as $h\to \infty$.

Applications of homological stability to arithmetic problems over function fields $F=\bF_q(B)$ go back to \cite{EVW}, and have been explored in a variety of contexts, e.g., Cohen-Lenstra heuristic, Malle's conjecture, and Manin's conjecture, see, e.g., \cite{EVW}, \cite{LL24a}, \cite{LL24b}, \cite{ETW}, \cite{LL25}, and \cite{DLTT25}.

It is natural to also consider {\em weak approximation}. 
In this context, weak approximation asserts the existence of sections matching a finite set of jet conditions, 
see \cite{GHS} for the existence of sections and \cite{HT-weak} for weak approximation; this translates into {\em existence} of $\bk$-points on spaces of such sections in the stable regime, when $h\to \infty$. Recall that {\em weak approximation} for $X$ over $F$ means that for any finite set 
of {\em admissible jets} there exists a section $\sigma: B\to \cX$ 
matching these jets;  
an admissible $N$-th jet is the truncation of a local analytic section of $\pi$ at $b\in B$ modulo $\mathfrak m_{B,b}^{N+1}$, power of the maximal ideal of the local ring at $b$. 
By \cite[Theorem 3]{HT-weak}, 
a geometrically rationally connected variety over $F $ satisfies {\em weak approximation} at places of good reduction; this has been extended to bad reduction places in some other situations, e.g.,  \cite{HT-ww}, \cite{xu-weak}, \cite{cubic}.   
{\em Effective} weak approximation, considered in \cite{HT12}, seeks effective control over the height of sections matching these jets. 
This is also coupled to {\em equidistribution}, which
can be viewed as a strong, quantitative form of weak approximation, and 
would require uniform control over the number of $\bk$-points on the corresponding spaces, for $\bk=\bF_q$, as in \eqref{eqn:sect}.  

\

Pursuing these analogies, we 
\begin{itemize} 
\item 
formulate
the relevant version of homological stability, and
\item 
prove it for fibrations arising from projectivizations of vector bundles, conic bundles, and quadric surface bundles over curves.
\end{itemize}

From now on, we assume that $\bk=\bC$ and identify varieties with their sets of complex points. Let 
$$
\alpha\in \rH_2(\cX(\bC), \bZ)
$$ 
be the class of a section of $\pi$; in particular, $\alpha\cdot \cX_b=1$, for all fibers $\cX_b=\pi^{-1}(b)$, $b\in B$. 
Let 
$$
\beta\in \rH_2(\cX(\bC),\bZ)
$$ 
be the class of a very free rational curve in a smooth fiber of $\pi$. 
Clearly, $\alpha+\beta$ is also the class of a section of $\pi$. 
Let 
$$
\Sigma=\{\hat{\sigma}_j\}_{j\in J}
$$ 
be a finite set of jets on $\cX$ in distinct fibers $\cX_{b_j}$, where each
$\hat{\sigma}_j$ is an admissible $N_j$-th jet. Let 
\begin{multline*} 
\Sect(\cX/B, \alpha,\Sigma):=
\\
\{\sigma :B\to \cX  \mid  \pi\circ \sigma = \mathrm{id}, \, \, [\sigma(B)]\equiv\alpha, \, \, \, \sigma \equiv \hat{\sigma}_j \pmod{\mathfrak m_{B,b_j}^{N_j+1}}, \, \,\forall j\in J
\} 
\end{multline*}
be the space of sections of class $\alpha$ matching all jets in $\Sigma$; this is a quasi-projective scheme over $\mathbb C$ and was constructed by Grothendieck (see, e.g., \cite[4(c)]{Gro95}). It can be realized as a closed subscheme of the morphism scheme $\mathrm{Mor}(B, \mathcal X, \alpha)$, which, in turn, is realized as a Zariski open subscheme of the Hilbert scheme of $\mathcal X$.

As explained in \cite{HT-weak}, an $N$-th jet condition in a fiber of $\cX$ over $b\in B$ can be rephrased in terms of a $0$-th jet condition on some component of a fiber over $b$ of an iterated blowup of $\cX$ with center in fibers over $b$. 
We know little about these spaces, e.g.,  
\begin{itemize}
    \item Are they irreducible and of expected dimension? 
    \item Is there some stabilization in their homology?
\end{itemize}
The first question has been studied in absence of jet conditions, e.g., in \cite{HRS}, \cite{RY16}, \cite{LT17}, \cite{LTJAG},\cite{LT19}, \cite{LT21}, \cite{BLRT20}, \cite{LRT23}, \cite{Okamura24}, \cite{Okamura25}, and \cite{BJ}. The second question has not been addressed, to our knowledge. This motivates the introduction of 

\ 

\noindent
{\bf Condition} {\bf (HS)} (Homological stability, for the triple $(\alpha, \beta, \Sigma)$): 
For $\alpha, \beta$, and $\Sigma$ as above, there
exists
a linear function $\ell$, with positive leading term, such that, for all $i\le \ell(m)$, one has 
$$
\rH_i(\Sect(\cX/B,\alpha+m\beta,\Sigma),\bZ)\simeq  \rH_{i}(\Sect(\cX/B,\alpha+(m + 1)\beta,\Sigma),\bZ). 
$$

\

Some geometric assumptions on $X$ will be necessary. E.g., could this hold for $X$ a del Pezzo surface? Or a toric variety? One of our main results is:

\begin{theo}
\label{thm:main}
    Homological stability, for any $(\alpha,\beta,\Sigma)$ with nonempty $\Sigma$,
    holds for:
    \begin{itemize}
\item $\bP(\cE)$, projectivization of a vector bundle $\cE$ of rank $\ge 2$ over $B$, 
\item smooth conic bundles over $B$,
\item smooth nonsplit quadric surface bundles over $B$, with at most $\mathsf A_1$-singular fibers. 
\end{itemize}
\end{theo}

Note that weak approximation and Manin's conjecture are known in many cases, e.g., for
quadric surfaces over $F=\bF_q(t)$. However, homological stability for spaces of sections (without jet matching conditions) has only been established for trivial families \cite{Segal}, \cite{Guest}, \cite{DT24}, \cite{DLTT25}.

\ 

In the same vein, when $\Sigma$ is a set of $0$-th jets, we consider the space 
$$
\Sect^{\rm top}(\cX/B,\alpha,\Sigma)
$$
of {\em topological sections} of corresponding classes; for $N_j$-th jets, with $N_j\ge 1$, we replace $\cX$ with a birational model realizing the jet as a point on some component of the fiber. 
This allows to formulate the parallel condition:

\ 

\noindent
{\bf Condition} {\bf (HST)} (Homological stability of topological sections):
The homology 
$$
\rH_i(\Sect^{\rm top}(\cX/B,\alpha+m\beta,\Sigma),\bZ)
$$
stabilizes, in the sense above. 

\

We establish Condition {\bf (HST)} in broad generality, see Corollary~\ref{coro:homologicalstabilitytop}. In applications to projective bundles, 
we show that the inclusion
\begin{equation}
\label{equation:algevstop}
\Sect(\cX/B,\alpha + m\beta,\Sigma) \hookrightarrow\Sect^{\rm top}(\cX/B,\alpha + m \beta,\Sigma)
\end{equation}
induces isomorphisms for low-degree homologies. This allows
to establish Condition {\bf (HS)} via Condition {\bf (HST)}, see   
Theorem~\ref{theo:projectivebundlealgvstop}.
The proofs rely on an explicit characterization of the relevant spaces of sections and topological gluing arguments, 
combined with tight control over the combinatorics of bar complexes, as in \cite{DT24, DLTT25}.

\bigskip

\noindent
{\bf Acknowledgments:}
The authors thank the referee for detailed comments which significantly improved the exposition.
The first author was partially supported by JST FOREST program Grant number JPMJFR212Z and by JSPS KAKENHI Grant-in-Aid (B) 23K25764. The second author was partially supported by 
NSF grant 2301983.

\section{Generalities}
\label{sect:gene}

\subsection*{Notation}

The cardinality of a finite set $J$ is denoted by $|J|$. 
For an admissible set of jets $\Sigma=\{ \hat{\sigma}_j\}_{j\in J}$, 
where $\sigma_j$ are $N_j$-th jets, 
we let
$$
\deg(\Sigma)=\sum_{j\in J}(N_j+1),  
$$
be its total degree.

Schemes are separated of finite type over $\bC$, and varieties are integral. 
For complex varieties, {\em dimension} refers to their complex dimension; for semi-algebraic sets, it is their real dimension. 
Isomorphisms of schemes or complex manifolds are denoted by $\cong$, homeomorphisms  by $\approx$, and homotopy equivalences by $\sim$.

\

In this section, we introduce a {\em semi-topological} model of the 
space of topological sections 
\[
\mathrm{Sect}^{\mathrm{top}}(\mathcal X/B, \alpha, \Sigma),
\]
based on ideas from \cite[Section 6]{DT24}. We make no assumptions on the geometry of the smooth projective variety $X$ over $\bC(B)$.

\subsection*{Compactly generated topologies}

Let $K$  be a compact space and $U\subset \bC$ an open subset containing $0$, with coordinate $z$. 
Let
$$
    \cC(K,U):=\{ f : K\times U\to \mathbb C\}
    $$ 
be the set of continuous functions, viewed as the set of continuous families, parametrized by $K$, 
of continuous functions $U\to \bC$.
This defines a compactly generated topology on $\cC(U)$, the set of continuous functions on $U$: 
a continuous map
$$
\bar{f} : K \to \cC(U)
$$
is equivalent to the condition that the corresponding family 
$$
f : K \times U \to \mathbb C
$$ 
satisfies $f\in \cC(K, U)$. 
Similar construction applies to
$$
\cC_0(K,U):=\{ f : K\times U \to \mathbb C \, \mid \, f(\xi,0)= 0, \quad  \forall \xi\in K\} \subset \cC(K,U).
$$
Fix a positive integer $N$ and a polynomial $\rho\in \bC[z]$ of degree $N$.
Let 
$$
 \cC(K,U,\rho,N) \subset\cC(K,U)
$$
be the set of continuous functions, congruent to $\rho(z)$ modulo $|z|^N$, i.e.,  for every $\xi \in K$, 
    \[
    f(\xi,z) - \rho(z) = o(|z|^N),\quad z \to 0,
    \]
    uniformly in $\xi \in K$.
As above, this induces a compactly generated topology on 
$$
\cC(U,\rho,N). 
$$

The following is analogous to \cite[Proposition 6.2]{DT24}:

\begin{lemm}
\label{lemm:semitopologicalcorrespondence}
    There is a natural bijection 
\begin{equation}
\label{eqn:bi}
\cC_0(K,U) \leftrightarrow \cC(K,U,\rho,N).
\end{equation}
\end{lemm}

\begin{proof}
    Given an $f_0\in \cC_0(K,U)$,  
    we have 
    \[
    f(\xi,z) = \rho(z) + z^Nf_0(\xi,z) \in \cC(K,U,\rho,N).
    \]
Conversely, for $f\in \cC(K,U,\rho,N)$, we put
    \[
    f_0(\xi,z) = 
    \begin{cases}
        (f(\xi,z) - \rho(z))/z^N & \text{if $z \neq 0$};\\
        0 & \text{if $z = 0$,}
    \end{cases}    
    \]
establishing the claim.
\end{proof}

A corollary of \eqref{eqn:bi} is the following:
\begin{coro}
\label{coro:homeo_functionspace}
There is a homeomorphism of topological spaces
\begin{equation} 
\label{eqn:bii}
\cC_0(U) \approx \cC(U,\rho,N). 
\end{equation}

\end{coro}

\begin{proof}
  Equation~\eqref{eqn:bi} establishes a bijection between  spaces.
  By properties of compactly generated topology, to verify that
  \[
\Phi : \cC_0(U) \to \cC(U,\rho,N)
  \]
  is continuous, it suffices to show that for any compact set $K$ and continuous map $\bar{f} : K \to \cC_0(U)$, the composition $$\Phi  \circ \bar{f} : K \to \cC_0(U, \rho, N)$$
  is also continuous.
  The map $\bar{f}$ induces  a continuous family of functions parametrized by $K$
  \[
  [f: K \times U \to \mathbb C] \in \cC_0(K,U).
  \]
  The map $\Phi  \circ \bar{f} : K \to \cC_0(U, \rho, N)$ corresponds to $\Phi_K(f)$, where 
  \[
   \Phi_K : \cC_0(K,U) \rightarrow \cC(K,U,\rho,N)
  \]
  is the map 
  established in Lemma~\ref{lemm:semitopologicalcorrespondence}, 
  and our assertion follows.
  \end{proof}

\subsection*{Semi-topological models}

The following is inspired by \cite[Definition 6.5]{DT24}:

\begin{itemize} 
\item For $b\in B$, let $U_b\subset B$ be an open 
neighborhood, with local holomorphic coordinate $z$ such that $z(b) = 0$. 
\item 
For a smooth $x\in \cX_b$, let $\mathcal U_x\subset \cX$ 
be an open neighborhood such that 
$\pi(\mathcal U_x) \subset U_b$.
\item 
Assume that we have local holomorphic coordinates $z, z_1, \ldots, z_n$ on $\mathcal U_x$ such that $\pi : \mathcal U_x \to U_b$ corresponds to 
$$
(z, z_1, \ldots, z_n) \mapsto z, \quad \text{ and } \quad z_1(x) = 0, \ldots, z_n(x) = 0.
$$
\item 
For $\hat{\sigma}$, an admissible $N$-th jet at $x$, let $\rho_i\in \bC[z]$ be a polynomial of degree $\leq N$ induced by $\hat{\sigma}$ and $z_i$, for $i=1,\ldots, n$. 
\end{itemize}
Let $K$ be a compact set and 
$$
\sigma : K \times B \to \mathcal X
$$ 
a family of continuous sections of $\pi$ parametrized by $K$. Since $K$ is compact, one can find an open $U_b' \subset U_b$ such that $\sigma(K \times U_b') \subset \mathcal U_x$. 

\begin{defi}
\label{defn:semi}
We say that $\sigma$ induces $\hat{\sigma}$ if for some (equivalently, every) 
$$
U_b, \, \mathcal U_x,\, z,\{ z_i\},\,  \{ \rho_i\}, \, U_b',
$$
one has 
\[
z_i\circ \sigma|_{K \times U'_b} \in \mathcal C(K, U'_b, \rho_i, N),\quad \text{for $i = 1, \ldots, n$.}
\]

Let $\Sigma = \{ \hat{\sigma}_j\}_{j\in J}$ be a set of admissible jets, where $\hat{\sigma}_j$ is an $N_j$-th jet.
Define
\[
\mathrm{Sect}^{\mathrm{stop}}(K,\mathcal X/B, \alpha, \Sigma)
\]
to be the set of continuous families, parametrized by $K$, of topological sections of $\pi$ of class $\alpha$
constrained by $\Sigma$, i.e., 
sections inducing $\hat{\sigma}_j$,  for all $j\in J$.
\end{defi}

Combining the argument of Proposition 6.9 and Lemma 6.10 of
\cite {DT24}, one verifies that this
defines a compactly generated topology on 
\[
\mathrm{Sect}^{\mathrm{stop}}(\mathcal X/B, \alpha, \Sigma).
\]
This is 
the {\em semi-topological model} of the quasi-projective scheme
$$
\mathrm{Sect}(\mathcal X/B, \alpha, \Sigma).
$$

\section{Properties of semi-topological models}
\label{sec:conditionHST}

\subsection*{Blowups}
\label{subsec:semitopvsblowup}

Consider a point $b \in B$ and an admissible $N$-th jet $\hat{\sigma}$, supported at a smooth point $x\in \mathcal X_b$, 
with $N \geq 1$. 
Let 
$$
\varphi: \tilde{\mathcal X}\to \cX
$$ 
be the blowup at $x$, with exceptional divisor $E$, and $\Sigma=\{\sigma\}$. 
Let  $\tilde{\alpha}$ be 
the class of a section of $\tilde{\cX}\to B$, with $\varphi_*\tilde{\alpha} = \alpha$ and $\tilde{\alpha}\cdot E = 1$.

\begin{lemm}
\label{lemm:blow}

The pushforward defines an isomorphism of algebraic varieties
    \[
    \varphi_* : \Sect(\tilde{\cX}, \tilde{\alpha}, \tilde{\Sigma}) \cong \Sect(\cX, \alpha, \Sigma),
    \]
    where $\tilde{\Sigma}$ consists of an admissible $(N-1)$-th jet on $\tilde{X}$. 
    Moreover, the pushforward induces a homeomorphism
    \[
    \varphi_*^{\rm{stop}} : \Sect^{\mathrm{stop}}(\tilde{\cX}, \tilde{\alpha}, \tilde{\Sigma}) \approx \Sect^{\mathrm{stop}}(\cX, \alpha, \Sigma).
    \]
\end{lemm}

\begin{proof}
    The first statement has been proved in \cite[Section 2.3]{HT-weak}. Let us prove the topological statement.
    As in Definition~\ref{defn:semi}, choose
    $$
U_b, \, \mathcal U_x,\, z,\{ z_i\},\,  U_b'.
    $$
Let $\rho_i\in \bC[z]$ be a polynomial of degree $\leq N$ induced by $\hat{\sigma}$ and $z_i$, and put 
$$
\rho_i'(z) = \rho_i(z)/z.
$$
Let $\sigma$ be a topological section of $\pi$ inducing $\hat{\sigma}$.
Then  
\[
z_i(\sigma(z)) -\rho_i(z) =  o(|z|^N), \quad i = 1, \ldots, n.
\]
The blowup $\tilde{\cX}$ is locally defined by
\[
\varphi^{-1}(\mathcal U_x) \cong Z(t_iz = sz_i, t_iz_j = t_jz_i) \subset \mathcal U_x \times \bP^n,
\quad i, j = 1, \ldots, n,
\]
where $(s : t_1 : \ldots : t_n)$ are homogeneous coordinates on $\bP^n$. Let $U'_b \subset U_b$ be an open neighborhood of $b$ such that $\sigma(U'_b) \subset \mathcal U_x$. We define the strict transform $\tilde{\sigma} : B \to \tilde{\cX}$ by 
\[
\tilde{\sigma}(b') = 
\begin{cases}
  \varphi^{-1}\circ \sigma(b') & \text{ if $b' \not\in U'_b$},\\
  (\sigma(b'))\times [z(b'): z_1(\sigma(b')): \ldots : z_n(\sigma(b'))] & \text{ if $b' \in U'_b \setminus \{b\}$},\\
  (\sigma(b)) \times [1: \rho_1'(0): \ldots: \rho_n'(0)]& \text{ if $b' = b$.}
\end{cases}
\]
Let $\tilde{x} = \tilde{\sigma}(b)$. Its local coordinates are $z, \tilde{z}_1, \ldots, \tilde{z}_n$, where $\tilde{z}_i = t_i/s$. We have
\[
\tilde{z}_i(\tilde{s}(z)) -\rho_i'(z)  =  o(|z|^{N-1}).
\]
Let $\hat{\sigma}^{(1)}$ be an admissible $(N-1)$-th jet of $\tilde{\cX}$ determined by $\rho_i'(z)$. 
The above shows that $\tilde{\sigma}$ induces $\hat{\sigma}^{(1)}$.
Conversely, if $\tilde{\sigma}$ induces $\hat{\sigma}^{(1)}$, then $\sigma = \varphi \circ \tilde{\sigma}$ satisfies 
\[
z_i(s(z)) = z\rho'_i(z)+ o(|z|^N) =  \rho_i(z) + o(|z|^N).
\]
This correspondence establishes that the map
\[
    \varphi_*^{\rm{stop}} : \Sect^{\mathrm{stop}}(\tilde{\cX}, \tilde{\alpha}, \tilde{\Sigma}) \to \Sect^{\mathrm{stop}}(\cX, \alpha, \Sigma).
    \]
is a bijection. 

To show that this is a homeomorphism, let $K$ be a compact set. The above discussion shows that the birational map $\varphi$ induces a map
\[
    \varphi_{K*}^{\rm{stop}} : \Sect^{\mathrm{stop}}(K, \tilde{\cX}, \tilde{\alpha}, \tilde{\Sigma}) \to \Sect^{\mathrm{stop}}(K, \cX, \alpha, \Sigma),
    \]
    and this is also bijective.
    Arguing as in the proof of Corollary~\ref{coro:homeo_functionspace}, we conclude that $\varphi_*^{\rm{stop}}$ is a homeomorphism.
\end{proof}

We call $\tilde{\Sigma}=\{ \hat{\sigma}^{(1)}\}$ the strict transform of $\Sigma=\{ \hat{\sigma}\}$ on $\tilde{\cX}$. Iterating this construction, for any set $\Sigma=\{ \hat{\sigma}_j\}_{j\in J}$ of admissible $N_j$-th jets, we produce an iterated blowup 
$$
\varphi_\Sigma: \tilde{\cX}_\Sigma\to \cX,
$$
with a section class $\tilde{\alpha}_\Sigma$ and 
a set $\tilde{\Sigma}$ of admissible $0$-th jets.

\begin{prop}
\label{prop:blowupmodel} 
    The pushforward map defines an isomorphism
    \[
    \varphi_{\Sigma*} : \Sect(\tilde{\cX}_\Sigma, \tilde{\alpha}_\Sigma, \tilde{\Sigma}) \cong \Sect(\cX, \alpha, \Sigma).
    \]
    Moreover, the pushforward map also induces a homeomorphism
    \[
    \varphi_{\Sigma*}^{\mathrm{stop}} : \Sect^{\mathrm{stop}}(\tilde{\cX}_\Sigma, \tilde{\alpha}_\Sigma, \tilde{\Sigma}) \approx \Sect^{\mathrm{stop}}(\cX, \alpha, \Sigma).
    \]
\end{prop}

\begin{proof}
    We apply Lemma~\ref{lemm:blow} iteratively, to successive blowups.
\end{proof}

\subsection*{Moving the support of the jet}

Let $\hat{\sigma}_0, \hat{\sigma}_0'$ be admissible $0$-th jets above $b_0$ which are supported in the {\em same} component of the fiber $\mathcal X_{b_0}$, but with disjoint support. 
Let $\Sigma_0=\{ \hat{\sigma}_j\}_{j\in J}$  be a set of admissible $N_j$-th jets over points $b_j$ distinct from $b_0$ and put
$$
\Sigma:=\{\hat{\sigma}_0\} \cup \Sigma_0, \quad \Sigma':=\{\hat{\sigma}_0'\} \cup \Sigma_0. 
$$

\begin{prop}
\label{prop:pointsindependnetgeneral}
We have a homeomorphism
\[
\Sect^{\mathrm{stop}}(\mathcal X, \alpha, \Sigma) \approx \Sect^{\mathrm{stop}}(\mathcal X, \alpha, \Sigma').
\]
\end{prop}

\begin{proof}
It suffices to show that the evaluation map
\[
\mathrm{ev} : \Sect^{\mathrm{stop}}(\mathcal X, \alpha, \Sigma_0) \to \mathcal X_{b_0}, \quad \sigma \mapsto \sigma(b_0)
\]
is topologically locally fiberwise isotrivial on the smooth locus of $\mathcal X_{b_0}$. 

Let $x_0,x_0'$ be the supports of $\hat{\sigma}_0, \hat{\sigma}_0'$.
Consider the following data:
\begin{itemize}
    \item $U_{b_0}$ an open neighborhood of $b_0$ such that $\overline{U}_{b_0}$ does not contain $b_j$, for $j\in J$;
    \item a homeomorphism of closures $\overline{U}_{b_0} \approx \overline{\mathbb D}$, where $\mathbb D$ is the open unit disk with coordinate $z$, mapping $U_{b_0}$ to $\mathbb D$ and $b_0$ to $0$;
    \item $\mathcal U_{x_0}$ an open neighborhood of $x_0$ such that $\pi(\mathcal U_{x_0}) = U_{b_0}$,
    \item a homeomorphism $\overline{\mathcal U}_{x_0} \approx   \overline{\mathbb D} \times \overline{\mathbb B}$, 
    where $\mathbb B$ is a unit ball with coordinates $w=(w_1,\ldots, w_n)$,
    yielding a commutative diagram

    \

\centerline{ 
    \xymatrix{ \overline{\cU}_{x_0} \ar[d]_{\pi}   \ar[r]^{\approx}   & \ar[d] \overline{\mathbb D} \times \overline{\mathbb B} \\ 
                        U_{b_0}      \ar[r]^{\approx}  &  \overline{\mathbb D} 
    }
}

\noindent 
where the right vertical map is projection onto the first factor.     
    \item we assume that $x_0,x_0'$ are contained in $\cX_{b_0}\cap \cU_{x_0}$,  $x_0=(0,0)$ in the coordinates on $\overline{\mathbb D} \times \overline{\mathbb B}$, 
    we suppose that $x_0'=(0,w')$. 
\end{itemize}
Our goal is to show that $\mathrm{ev}^{-1}(x_0)\approx \mathrm{ev}^{-1}(x_0')$. Let 
$$
\rho_0 : \overline{\mathbb B} \to \overline{\mathbb B}
$$ 
be a homeomorphism  inducing the identity on the boundary and 
mapping $0$ to $w'$. 
We construct a homotopy
$$
\rho_t : [0,1] \times \overline{\mathbb B} \to  \overline{\mathbb B}$$ such that 
\begin{itemize}
    \item $\rho_1$ is the identity and $\rho_0$ is the map above, and
    \item for any $t \in [0, 1]$, $\rho_t : \overline{\mathbb B} \to \overline{\mathbb B}$ is a homeomorphism
    inducing identity on the boundary.
\end{itemize}
Then we define a $\pi$-homeomorphism $\Phi : \mathcal X \to \mathcal X$ by
\[
\Phi(x) = 
\begin{cases}
    x & \text{if $x \not\in \overline{\mathcal U}_{x_0}$}, \\
    (z(\pi(x)), \rho_{|z(\pi(x)|}(w(x))) & \text{if $x \in \overline{\mathcal U}_{x_0}$}. 
\end{cases}
\]
This realizes $\mathrm{ev}^{-1}(x_0)\approx \mathrm{ev}^{-1}(x_0')$, as it maps $x_0$ to $x_0'$. 
Indeed, we obtain a bijection
\[
\Sect^{\mathrm{stop}}(\mathcal X, \alpha, \Sigma) \to  \Sect^{\mathrm{stop}}(\mathcal X, \alpha, \Sigma'), [\sigma: B \to \mathcal X] \mapsto \Phi\circ\sigma.
\]
Similarly, for any compact set $K$, the map
\[
\Sect^{\mathrm{stop}}(K, \mathcal X, \alpha, \Sigma) \to  \Sect^{\mathrm{stop}}(K, \mathcal X, \alpha, \Sigma'), [\sigma: K\times B \to \mathcal X] \mapsto \Phi\circ\sigma,
\]
is a bijection, and we conclude as in Corollary~\ref{coro:homeo_functionspace}.
\end{proof}

\subsection*{Gluing a rational curve}

Now we assume that the fiber $\cX_{b_0}$ is smooth and that it contains a rational curve of class $\beta$
joining $x_0,x_0'$.

\begin{prop}
\label{prop:stabilityHST}
We have a homotopy equivalence:
 \[
 \Sect^{\mathrm{stop}}(\mathcal X/B, \alpha,\Sigma) \sim \Sect^{\mathrm{stop}}(\mathcal X/B, \alpha + \beta, \Sigma).
 \]
\end{prop}

\begin{proof} 

Let $f : \mathbb P^1 \to \mathcal X_{b_0}$ be a rational curve of class $\beta$, such that 
$f([1:0]) = x_0$ and $f([0:1]) = x_0'$. We introduce continuous functions:
$$
\eta :  [0,1] \to [0,1],
\quad 
\eta(t) := 
\begin{cases}
    0 & \text{ if $t \in [0, \frac{1}{2}]$}, \\
    2(t-\frac{1}{2}) & \text{ if $t \in [\frac{1}{2}, 1]$,}
\end{cases}
$$
and 
\[
\zeta : [0, 1/2) \to [0,+ \infty), \quad  \zeta(t):=\tan(\pi t).
\]
    Let $U_{b_0}$ be an open neighborhood of $b_0$, with a fixed homeomorphism to $\bD$, as above, with coordinate $z$ and center $b_0$. Assume that   
    \begin{itemize}
    \item $\overline{U}_{b_0}$ does not contain $b_j$, for all $j\in J$,
    \item 
    there is a $\overline{\bD}$-homeomorphism 
    $$
    \phi : \mathcal X|_{\overline{U}_{b_0}} \approx \overline{\bD} \times \mathcal X_{b_0}
    $$ 
    such that $\phi(x_0) = (0, x_0)$.
    \end{itemize}
    Let $\sigma \in \mathrm{Sect}^{\mathrm{stop}}(\cX/B, \alpha, \Sigma)$. We define 
    $
    \sigma_1 \in \mathrm{Sect}^{\mathrm{stop}}(\cX/B, \alpha,\Sigma)
    $
    by 
    \begin{equation}
    \sigma_1(b) := 
    \begin{cases}
        \sigma(b) & \text{ if $b \not\in \overline{\bD}$,}\\
        \phi^{-1}\left(z(b), \phi_2\left(\sigma\left(\frac{\eta(|z(b)|)}{|z(b)|}z(b)\right)\right)\right) & \text{ if $b \in \overline{\bD}$ and $b \neq b_0$,}\\
        x_0 & \text{ if $b = b_0$,}
    \end{cases}
    \label{equation:sigma_1}
    \end{equation}
     where $\phi_2$ is the composition of $\phi$ with the second projection.
     When $|z(b)|\leq 1/2$, we have
     $$\phi_2(\sigma_1(b)) = x_0.$$
    Next we define $\sigma' \in \mathrm{Sect}^{\mathrm{stop}}(\cX/B, \alpha + \beta, \Sigma')$ by 
    \[
    \sigma'(b)\!\! := \!\!
    \begin{cases}
        \sigma_1(b) & \!\!\text{if $b \not\in \overline{\bD}$,}\\
        \sigma_1(b) & \!\!\text{if $b\in \overline{\bD}$ and $|z(b)| \geq 1/2$,}\\
        \phi^{-1}(z(b), f([\zeta(|z(b)|)z : 1])& \!\!\text{if $b\in \overline{\bD}$ and $|z(b)| < 1/2$.}
    \end{cases}
    \]
    One may think of $\sigma'$ as obtained from $\sigma$ by gluing the vertical rational curve
    $$
   f :  \mathbb P^1  \to \mathcal X_{b_0}
    $$ 
    at $x_0$. By construction,  $\sigma'(b_0) = x_0'$.
    This defines a continuous map
    \[
    \Phi : \mathrm{Sect}^{\mathrm{stop}}(\cX/B, \alpha, \Sigma) \to 
    \mathrm{Sect}^{\mathrm{stop}}(\cX/B, \alpha + \beta, \Sigma'),\quad  \sigma \mapsto \sigma'.
    \]
    Indeed, this is continuous because the same construction establishes
    \[
    \Phi_K : \mathrm{Sect}^{\mathrm{stop}}(K, \cX/B, \alpha, \Sigma) \to 
    \mathrm{Sect}^{\mathrm{stop}}(K, \cX/B, \alpha + \beta, \Sigma'),\quad  \sigma \mapsto \sigma',
    \]
    for any compact $K$.
    
We construct the homotopy inverse to $\Phi$ by gluing an inversely oriented sphere so that the composition is homotopic to the identity.
Indeed, for $\tau' \in \mathrm{Sect}^{\mathrm{stop}}(\cX/B, \alpha + \beta, \Sigma')$, we construct $\tau_1' \in \mathrm{Sect}^{\mathrm{top}}(\cX/B, \alpha + \beta, \Sigma')$ in the same way we constructed $\sigma_1$. In particular, when $|z(b)|\leq 1/2$, we have
     $$\phi_2(\tau_1'(b)) = x_0'.$$
     We define $\tau\in \mathrm{Sect}^{\mathrm{top}}(\cX/B, \alpha,   \Sigma)$ by 
    \[
    \tau(b) = 
    \begin{cases}
        \tau_1'(b) & \text{$b \not\in \overline{\bD}$,}\\
        \tau_1'(b) & \text{$b\in \overline{\bD}$ and $|z(b)| \geq 1/2$,}\\
        \phi^{-1}(z(b), f([1: \zeta(|z(b)|)\overline{z}])& \text{$b\in \overline{\bD}$ and $|z(b)| < 1/2$.}
    \end{cases}
    \]
   Complex conjugation corresponds to gluing the inversely oriented sphere. This defines a continuous map
    \[
    \Psi : \mathrm{Sect}^{\mathrm{stop}}(\cX/B, \alpha+\beta, \Sigma') \to \mathrm{Sect}^{\mathrm{stop}}(\cX/B, \alpha, \Sigma), \quad \tau' \mapsto \tau.
    \]
    Since gluing of a sphere and the inversely oriented sphere is homotopic to a point, the compositions
    \[
    \Psi \circ \Phi, \quad \Phi \circ \Psi,
    \]
    are homotopic to identities.
    Indeed, pick $\sigma \in \mathrm{Sect}^{\mathrm{stop}}(\cX/B, \alpha, \Sigma)$. Then we have
     \[
    \Phi(\sigma)(b)\!\! := \!\!
    \begin{cases}
        \sigma(b) & \!\!\text{if $b \not\in \overline{\bD}$,}\\
        \phi^{-1}\left(z(b), \phi_2\left(\sigma\left(\frac{\eta(|z(b)|)}{|z(b)|}z(b)\right)\right)\right) & \!\!\text{if $b\in \overline{\bD}$ and $|z(b)| \geq 1/2$,}\\
        \phi^{-1}(z(b), f([\zeta(|z(b)|)z : 1])& \!\!\text{if $b\in \overline{\bD}$ and $|z(b)| < 1/2$.}
    \end{cases}
    \]
    Hence $\Psi\circ\Phi(\sigma)(b)$ is given by
    \[
    \begin{cases}
        \sigma(b) & \!\!\text{if $b \not\in \overline{\bD}$,}\\
        \phi^{-1}\left(z(b), \phi_2\left(\sigma\left(\eta(\eta(|z(b)|))z(b)/|z(b)|\right)\right)\right) & \!\!\text{if $b\in \overline{\bD}$, $|z(b)| \geq 3/4$,}\\
        \phi^{-1}(z(b), f([\zeta(\eta(|z(b)|))\frac{\eta(|z(b)|)}{|z(b)|}z(b) : 1])& \!\!\text{if $b\in \overline{\bD}$, $1/2\leq |z(b)| \leq 3/4$,}\\
         \phi^{-1}(z(b), f([1: \zeta(|z(b)|)\overline{z}])& \!\!\text{if $b\in \overline{\bD}$ and $|z(b)| < 1/2$.}
    \end{cases}
    \]
    We denote $\zeta(\eta(t))\eta(t)t$ by $\delta(t)$.
    Then we have
    \[
    \left[\zeta(\eta(|z|))\frac{\eta(|z|)}{|z|}z : 1\right] = [\delta(|z|) : \bar{z}] = \left[1:\frac{1}{\delta(|z|)} \bar{z}\right].
    \]
    For any $t \in [0,1]$, we define a homotopy
    \[
    F : [0,1] \times \mathrm{Sect}^{\mathrm{stop}}(\cX/B, \alpha, \Sigma) \to \mathrm{Sect}^{\mathrm{stop}}(\cX/B, \alpha, \Sigma), 
    \]
    where $F_t(\sigma)$ is given by
    \[
    \begin{cases}
        \sigma(b) & \!\!\text{if $b \not\in \overline{\bD}$,}\\
        \phi^{-1}\left(z(b), \phi_2\left(\sigma\left(\eta(\eta(|z(b)|))z(b)\right)\right)\right) & \!\!\text{if $b\in \overline{\bD}$,  $|z(b)| \geq 3/4$,}\\
        \phi^{-1}(z(b), f([1:\min\{\frac{1}{\delta(|z|)}, \zeta(t/2)\}\bar{z}])& \!\!\text{if $b\in \overline{\bD}$, $1/2\leq |z(b)| \leq 3/4$,}\\
         \phi^{-1}(z(b), f([1: \min\{\zeta(|z(b)|), \zeta(t/2)\}\overline{z}])& \!\!\text{if $b\in \overline{\bD}$,  $|z(b)| < 1/2$.}
    \end{cases}
    \]
    This is continuous because for any compact set $K$ and a continuous map 
    $$
    \bar{f} : K \to [0,1] \times \mathrm{Sect}^{\mathrm{stop}}(\cX/B, \alpha, \Sigma),
    $$
    the above construction maps the corresponding family
    \[
    f : [0,1]\times K \times B \to \mathcal X
    \]
    to a continuous family
    \[
     F(f) : [0,1]\times K \times B \to \mathcal X.
    \]
    Then we conclude that the corresponding map 
    \[
    \overline{F(f)} = F\circ\bar{f} : K \to [0,1] \times  \mathrm{Sect}^{\mathrm{stop}}(\cX/B, \alpha, \Sigma)
    \]
    is continuous.
We have $F_1 = \Psi \circ \Phi$ and
     \[
    F_0(\sigma) = 
    \begin{cases}
        \sigma(b) & \!\!\text{if $b \not\in \overline{\bD}$,}\\
        \phi^{-1}\left(z(b), \phi_2\left(\sigma\left(\eta(\eta(|z(b)|))z(b)\right)\right)\right) & \!\!\text{if $b\in \overline{\bD}$,  $|z(b)| \geq 3/4$,}\\
        \phi^{-1}(z(b), f([1:0]))& \!\!\text{if $b\in \overline{\bD}$,  $|z(b)| \leq 3/4$.}
    \end{cases}
    \]
    This $F_0$ is homotopy equivalent to the identity.  
    We conclude that 
    \[
    \Phi : \mathrm{Sect}^{\mathrm{stop}}(\cX/B, \alpha, \Sigma) \sim
    \mathrm{Sect}^{\mathrm{stop}}(\cX/B, \alpha + \beta,  \Sigma' )
    \]
    is a homotopy equivalence.
        Proposition~\ref{prop:pointsindependnetgeneral} implies that
         \[
        \mathrm{Sect}^{\mathrm{stop}}(\cX/B, \alpha, \Sigma) \sim
    \mathrm{Sect}^{\mathrm{stop}}(\cX/B, \alpha+\beta, \Sigma).
        \]
    \end{proof}

\subsection*{Independence of jets in smooth fibers}

Now we assume that the fiber $\cX_0$ is smooth, that 
$\hat{\sigma}_0$ and $\hat{\sigma}_0'$ are arbitrary admissible jets such that
the supports $x_0,x_0'$ of $\hat \sigma_0$ and $\hat \sigma_0'$ coincide, but the jets are distinct.

\begin{prop}
\label{prop:independentofjets}
We have a homotopy equivalence
\[
\Sect^{\mathrm{stop}}(\mathcal X/B, \alpha,  \Sigma) \sim \Sect^{\mathrm{stop}}(\mathcal X/B, \alpha ,  \Sigma').
\]
\end{prop}

We start with the following:
\begin{lemm}
\label{lemm:onlydependinglength}
    In the setting of Proposition~\ref{prop:independentofjets}, assume that $\hat{\sigma}_0$ and $\hat{\sigma}_0'$ have the same length.
    Then we have a homeomorphism
    \[
\Sect^{\mathrm{stop}}(\mathcal X/B, \alpha, \Sigma) \approx \Sect^{\mathrm{stop}}(\mathcal X/B, \alpha , \Sigma').
\]
\end{lemm}

\begin{proof}
    Let $\tilde{\mathcal X}_{\Sigma} \to B$ be the blowup model associated to $(\mathcal X/B, \Sigma)$ as in Proposition~\ref{prop:blowupmodel}. Let $\tilde{\Sigma}$ be the strict transform of $\Sigma$ which is the set of admissible $0$-th jets. Similarly we construct $\tilde{\mathcal X}_{\Sigma'} \to B$ and $\tilde{\Sigma}'$ associated to $(\mathcal X/B, \Sigma')$.
    By Proposition~\ref{prop:blowupmodel}, we have
    \[
    \Sect^{\mathrm{stop}}(\mathcal X/B, \alpha, \Sigma) \approx \Sect^{\mathrm{stop}}(\tilde{\mathcal X}_{\Sigma}/B, \alpha, \tilde{\Sigma}),
    \]
    and 
    \[
\Sect^{\mathrm{stop}}(\mathcal X/B, \alpha, \Sigma') \approx \Sect^{\mathrm{stop}}(\tilde{\mathcal X}_{\Sigma'}/B, \alpha, \tilde{\Sigma'}).
\]
Thus it suffices to show that
\[
\Sect^{\mathrm{stop}}(\tilde{\mathcal X}_{\Sigma}/B, \alpha, \tilde{\Sigma})  \approx \Sect^{\mathrm{stop}}(\tilde{\mathcal X}_{\Sigma'}/B, \alpha, \tilde{\Sigma'}).
\]
    Note that we can construct a proper smooth algebraic deformation over $B$ from $\tilde{\mathcal X}_{\Sigma}$ to $\tilde{\mathcal X}_{\Sigma'}$, in particular, this shows that $\tilde{\mathcal X}_{\Sigma}$ is $B$-homeomorphic to $\tilde{\mathcal X}_{\Sigma'}$. 
    Our assertion follows from Proposition~\ref{prop:pointsindependnetgeneral}.
\end{proof}

\begin{proof}[Proof of Proposition~\ref{prop:independentofjets}]
    It suffices to show the claim when $\hat{\sigma}_0$ is a $0$-th jet and $\hat{\sigma}_0'$ is an $N$-th jet. We assume that both are supported at $x_0 \in \mathcal X_{b_0}$.
    Let $\bD$ be an open neighborhood of $b_0$ such that
    \begin{itemize}
    \item $\overline{\bD}$ does not contain any $b_j, j\in J$,
    \item the closure $\overline{\bD}$ is homeomorphic to the closed unit disk with  complex coordinate $z$ and center corresponding to $b_0$, and;
    \item 
    there is a $\overline{\bD}$-homeomorphism $$
    \phi: \mathcal X|_{\overline{\bD}} \approx \overline{\bD} \times \mathcal X_{b_0}
    $$ 
    such that $\phi(x_0) = (0, x_0)$.
    By \cite[Proposition 9.5]{Voi02}, we may assume that the fibers of the composition $\phi_2 = \mathrm{pr}_2 \circ \phi$ are complex manifolds.
    \end{itemize}
     Let $\mathcal U_x \subset \mathcal X$ be an open neighborhood of $x$
     with local holomorphic coordinates $z, z_1, \ldots, z_n$ such that $\pi : \mathcal U_x \to \bD$ corresponds to mapping 
     \[
     (z, z_1, \ldots, z_n) \mapsto z \text{ and } z_1(x) = 0, \ldots, z_n(x) = 0.
     \]
     We may assume that the section 
     \[
     z \mapsto (z, 0, \ldots, 0)
     \]
     corresponds to the section 
     \[
     z \mapsto \phi^{-1}(z, x_0).
     \]
     This is possible because the fibers of $\phi_2$ are complex curves.
     
     Let $\rho_i\in \bC[z]$ be a polynomial of degree $\leq N$ induced by $\hat{\sigma}_0'$ and $x_i$. By Lemma~\ref{lemm:onlydependinglength}, we may assume that $\rho_i(z) = 0$, for all $i = 1, \ldots, n$.

     It is clear that 
     \[
     \Sect^{\mathrm{stop}}(\mathcal X/B, \alpha , \Sigma') \subset \Sect^{\mathrm{stop}}(\mathcal X/B, \alpha, \Sigma),
     \]
     so we may define 
     \[
     \Phi : \Sect^{\mathrm{stop}}(\mathcal X/B, \alpha , \Sigma') \to \Sect^{\mathrm{stop}}(\mathcal X/B, \alpha, \Sigma), \quad \sigma \mapsto \sigma.
     \]
     We construct the homotopy inverse. Let 
     $$
     \sigma  \in \Sect^{\mathrm{stop}}(\mathcal X/B, \alpha, \Sigma).$$ 
     Let $\sigma_1 : B \to \mathcal X$ be the section as constructed in (\ref{equation:sigma_1}). By  construction, $\sigma_1 \in \Sect^{\mathrm{stop}}(\mathcal X/B, \alpha , \Sigma')$. The inverse map is defined by
     \[
     \Psi : \Sect^{\mathrm{stop}}(\mathcal X/B, \alpha, \Sigma) \to \Sect^{\mathrm{stop}}(\mathcal X/B, \alpha , \Sigma'), \quad \sigma \mapsto \sigma_1. 
     \]
     It is easy to show that $\Phi \circ \Psi$ and $\Psi \circ \Phi$ are homotopic to identities. Indeed, such a homotopy can be obtained using the function
     \[
\eta_t(s) :  (0,1] \times [0,1] \to [0,1]: (s, t) \mapsto \eta_t(s),
\]
defined by
\[
\eta_t(s) := 
\begin{cases}
    0 & \text{ if $s \in (0, \frac{t}{2}]$}, \\
    (1-\frac{t}{2})^{-1}(s-\frac{t}{2}) & \text{ if $s \in [\frac{t}{2}, 1]$.}
\end{cases}
\]
\end{proof}

Combining Propositions~\ref{prop:pointsindependnetgeneral}, \ref{prop:stabilityHST}, and \ref{prop:independentofjets}, we obtain:

\begin{coro}
\label{coro:homologicalstabilitytop}
Let $\Sigma$ be an admissible jet datum with at least one jet supported in a smooth fiber $\mathcal X_{b_0}$.
Let $\beta$ be the class of a rational curve in $\mathcal X_{b_0}$.
Then we have a homotopy equivalence
\[
\Sect^{\mathrm{stop}}(\mathcal X/B, \alpha, \Sigma) \sim \Sect^{\mathrm{stop}}(\mathcal X/B, \alpha + \beta , \Sigma).
\]
\end{coro}

\section{Projective bundles over curves and the Abel--Jacobi map}
\label{sect:proj}

Let $B$ be a smooth projective irreducible curve of genus $g(B)$
and 
$$
\pi : \mathcal X = \mathbb P(\mathcal E) \to B
$$ 
the projectivization of a vector bundle $\cE$ over $B$ of rank $(n + 1)\ge 2$.
A class $\alpha$ of sections is specified by the degree 
$$
d = d(\alpha):=[\mathcal O_{\mathbb P(\mathcal E)/B}]\cdot \alpha;
$$
we denote the corresponding space of sections by
\[
\mathrm{Sect}(\mathcal X/B, d),
\]
and the subspace of sections with prescribed admissible jet data $\Sigma$ by
\[
\mathrm{Sect}(\mathcal X/B, d, \Sigma)\subseteq \mathrm{Sect}(\mathcal X/B, d).
\]

\subsection*{Abel--Jacobi map}
There is a well-known bijection 
$$
\{\sigma : B \to \cX\} \quad \leftrightarrow \quad \{ \mathcal E \twoheadrightarrow L\},  
$$
induced by
$$
\sigma \mapsto L:=\sigma^*\mathcal O_{\mathbb P(\mathcal E)/B}(1). 
$$
Consider the projections

\centerline{
\xymatrix{
& \ar[dl]_{\varpi} B\times \Pic^d(B) \ar[dr]^{\varpi_d} & \\
B & & \Pic^d(B) 
}
}

\noindent 
and let 
$$
\mathcal L_d\to B \times \Pic^d(B)
$$ 
be a universal line bundle. 
When $d\gg 0$ (depending on $\mathcal E$), the sheaf
$$
\mathcal V_d := \varpi_{d *}(\varpi^*(\mathcal E^\vee) \otimes \mathcal L_d)
$$ 
is locally free of rank 
$$
r = d(n + 1) - \deg(\mathcal E) + (n + 1)(1-g(B)).
$$

We have 

\ 

\centerline{
\xymatrix{
\mathrm{Sect}(\cX/B, d) \ar@{^{(}->}[r] \ar[dr]_{\mathrm{AJ}} & \bP(\cV_d^\vee) \ar[d] 
\\
                   & \Pic^d(B)
                   }
}

\noindent 
as a Zariski open subset of a projective bundle over $\Pic^d(B)$, where 
$\mathrm{AJ}$ is 
the Abel--Jacobi map.
Let 
$$
\tilde{S}_d \subset \cV_d
$$ 
be the Zariski open subset consisting of those $s \in \rH^0(B, \mathcal E^\vee \otimes L_d)$ which do not vanish on $B$; here $L_d\in\Pic^d(B)$ is the restriction of
$\cL_d$ to $B \times\{L_d\} \subset B \times \Pic^d(B)$. 
This is a $\mathbb C^\times$-torsor
\[
\tilde{S}_d \to \Sect(\cX/B,d).
\]

We turn to $\mathrm{Sect}(\mathcal X/B, d, \Sigma)$. 

\begin{lemm}
A section $\sigma : B \to \cX$ with $\deg(L)=d$
inducing $\Sigma$  corresponds to a section 
$$
s \in \rH^0(B, \mathcal E^\vee \otimes L_{d})
$$
satisfying the jet condition imposed by $\Sigma$, which is 
a linear condition.
Moreover, for  $d\gg 0$, depending on 
$\mathcal E$ and  $\deg (\Sigma)$, 
such sections form a codimension $n\cdot\deg (\Sigma)$ vector subbundle 
$$
\mathcal V_{d, \Sigma} \subset \mathcal V_d.
$$
\end{lemm}
\begin{proof}
Put $y_j= (N_j + 1)b_j$, an effective divisor on $B$. The homomorphism
\[
\mathcal E^\vee \to \oplus_{j \in J} (\mathcal E^\vee|_{y_j}),
\]
is surjective. The jet $\hat{\sigma}_j$ determines a codimension $(N_j + 1)n$ subspace $\mathcal H_j$ of $\mathcal E^\vee|_{y_j}$, and 
we denote the kernel of
\[
\mathcal E^\vee \to \oplus_{j \in J} (\mathcal E^\vee|_{y_j}/\mathcal H_j)
\]
by $\mathcal E^\vee(-\Sigma)$. 
We identify 
$$
(\mathcal V_{d, \Sigma})_{L_d}:=\rH^0(B, \mathcal E^\vee(-\Sigma)\otimes L_d),
$$ 
which is the kernel of
\[
\rH^0(B, \mathcal E^\vee(-\Sigma)\otimes L_d) \to \oplus_{j \in J} (\mathcal E^\vee|_{y_j}/\mathcal H_j).
\]
This shows that the incidence conditions imposed by the jet $\Sigma$ are linear conditions.

For  $d\gg 0$,  depending on $\mathcal E$ and $\deg(\Sigma)$, 
\[
\rH^0(B, \mathcal E^\vee\otimes L_d) \to \oplus_{j \in J} (\mathcal E^\vee|_{y_j}),
\]
is surjective. We conclude that $(\mathcal V_{d, \Sigma})_{L_d}$ has dimension independent of $L_d$, in this range. 
Thus we define 
\[
\mathcal V_{d, \Sigma} = (\varpi_d)_*(\varpi^*(\mathcal E^{\vee}(-\Sigma)) \otimes \mathcal L_d).
\]
\end{proof}

We have

\ 

\centerline{
\xymatrix{
\mathrm{Sect}(\cX/B,d,\Sigma) \ar@{^{(}->}[r] \ar[dr]_{\mathrm{AJ}} & \mathbb P(\mathcal V_{d, \Sigma}^\vee)  \ar[d] 
\\
                   & \Pic^d(B).
                   }
}

\ 

\noindent 
as a Zariski open subset. 
Similarly, we consider the Zariski open subset $\tilde{S}_{d, \Sigma} \subset \mathcal V_{d, \Sigma}$, yielding $\mathbb C^\times$-torsor
\[
\tilde{S}_{d, \Sigma} \to \mathrm{Sect}(\cX/B, d, \Sigma).
\]

\begin{prop}
\label{prop:projectivebundle}
For $d\gg 0$, depending on $\mathcal E$ and $\deg(\Sigma)$,  the space  
    \[
\mathrm{Sect}(\mathcal X/B, d, \Sigma)
\]
is a Zariski open subset of a projective bundle over $\Pic^d(B)$, of relative dimension
\[
d(n + 1) - \deg(\mathcal E) + (n + 1)(1-g(B)) -\sum_{j\in J} n(N_j+1) -1.
\]
In particular, it is irreducible, of expected dimension
\[
d(n + 1) - \deg(\mathcal E) + n(1-g(B)) -\sum_{j\in J} n(N_j+1).
\]
\end{prop}

\subsection*{Semi-topological counterparts}

Put
\[
\mathcal V^{\mathrm{stop}}_d := \{ (s, L) \, | \, L \in \Pic^d(B), s \in \rH^0_{\mathrm{cont}}(B, \mathcal E^\vee \otimes L) \},
\]
where $\rH^0_{\mathrm{cont}}$ denotes the space of continuous sections, with its
compactly generated topology.
Note that $\mathcal V^{\mathrm{stop}}_d$ is a locally trivial bundle of Banach spaces over $\Pic^d(B)$.
Let 
$$
\tilde{S}^{\mathrm{stop}}_d \subset \mathcal V^{\mathrm{stop}}_d
$$ 
be the open subset of $(s, L)$ such that $s$ is nowhere vanishing on $B$; it carries a $\bC^\times$-action
via 
\begin{equation}
    \label{eqn:c}
\lambda \cdot (s, L) = (\lambda s, L).
\end{equation} 
We have an inclusion
\[
\mathrm{Sect}(\cX/B, d) \hookrightarrow S^{\mathrm{stop}}_d:=\tilde{S}^{\mathrm{stop}}_d/\bC^\times. 
\]

We perform the same for jet conditions.
Fix $\Sigma =\{\hat{\sigma}\}_{j\in J}$ such that $\hat{\sigma}_j$ is supported at  $x_j \in \mathcal X_{b_j}= \mathbb P(\mathcal E_{b_j})$,  the fiber at $b_j \in B$. Let $\epsilon_{b_j} \subset \mathcal E_{b_j}^\vee$ be the $1$-dimensional subspace corresponding to $x_j$. 
Choose a Euclidean open neighborhood $U \subset \Pic^d(B)$
and a finite Euclidean open covering $\{U_\lambda\}$ of $B$ with  holomorphic trivializations
\begin{equation}
\label{equation:holomorphictrivialization}
\mathcal E|_{U_{\lambda}} \cong \oplus_{i = 1}^{n+1}\mathcal O, \quad \mathcal L_d|_{U \times U_{\lambda}} \cong \mathcal O.
\end{equation}
We construct a Banach subbundle
\[
\mathcal V^{\mathrm{stop}}_{d, \Sigma} \subset \mathcal V^{\mathrm{stop}}_{d},
\]
parametrizing continuous sections matching jets by specifying  the fiber  
over $L\in \Pic^d(B)$. 
Suppose we have a continuous section 
$$
s \in \rH^0_{\mathrm{cont}}(B, \mathcal E^\vee \otimes L).
$$
For each $b_j$, pick $\lambda$ such that $b_j \in U_{\lambda}$.
After shrinking $U_{\lambda}$ if necessary, we pick a homogeneous coordinate $z$ such that $z(b_j) = 0$ and $z$ exhibits a holomorphic isomorphism $U_{\lambda} \cong \mathbb D$ to the unit disk.

We choose the trivialization in (\ref{equation:holomorphictrivialization}) so that $\epsilon_{b_j}$ corresponds to 
\[
\mathbb C \cdot (1, 0, \ldots, 0) \subset \oplus_{i = 1}^{n+1}\mathcal O.
\]
The trivialization (\ref{equation:holomorphictrivialization}) induces local holomorphic coordinates 
\[
z, t_0, \cdots, t_n,
\]
of the bundle $\mathcal E^\vee|_{U_\lambda}$,
and local coordinates 
\[
z, z_1 = t_1/t_0, \ldots, z_n = t_n/t_0,
\]
in a neighborhood of $x_j$ in $\mathbb P(\mathcal E|_{U_{\lambda}})$.
The jet data define $\rho_i^j\in \bC[z]$.
We define a Banach bundle 
\[
\mathcal V^{\mathrm{stop}}_{d, \Sigma} \to \Pic^d(B)
\]
so that the fiber 
\[
(\mathcal V^{\mathrm{stop}}_{d, \Sigma})_L, \quad L \in \Pic^d(B),
\]
is the space of sections $s \in \rH^0_{\mathrm{cont}}(B, \mathcal E^\vee \otimes L)$ such that for each $j$, $s$ induces continuous functions $s^j_0, \ldots, s^j_n$, with respect to the trivialization (\ref{equation:holomorphictrivialization}), such that 
\[
\rho_i^j(z)s^j_0(z) = s^j_i(z)  + o(|z|^N), \quad i = 1, \ldots, n.
\]
These are linear conditions in $$s \in \rH^0_{\mathrm{cont}}(B, \mathcal E^\vee \otimes L).$$

Proof of Lemma~\ref{lemm:semitopologicalcorrespondence} shows that the  ${s}^j_i$ correspond to functions ${\check{s}^j_i}$ such that ${\check{s}^j_i}(0) = 0$, for $i = 1, \ldots, n$, via 
\[
\begin{pmatrix}
    {\check{s}^j_0}\\ {\check{s}^j_1} \\ \vdots \\ \check{s}^j_n
\end{pmatrix}
= 
\begin{pmatrix}
    1 & 0 & \cdots & 0\\
    -\frac{\rho^j_1}{z^{N_j}} & \frac{1}{z^{N_j}} & \cdots & 0\\
    \vdots & \vdots & \vdots & \vdots \\
    -\frac{\rho^j_n}{z^{N_j}} & 0 & \cdots & \frac{1}{z^{N_j}} 
\end{pmatrix}
\begin{pmatrix}
    s^j_0\\ s^j_1 \\ \vdots \\ s^j_n
\end{pmatrix}.
\]
Using this matrix as a transition matrix, we obtain the twisted holomorphic vector bundle $\check{\mathcal E}^\vee(-\Sigma)$ 
over $\Pic^d(B)$ such that 
$(\mathcal V^{\mathrm{stop}}_{d, \Sigma})_L$ is identified with 
\begin{multline*}
\rH^0_{\mathrm{cont}}(B, \check{\mathcal E}^\vee(-\Sigma) \otimes L)_\Sigma:= \\
\{ \check{s}\in \rH^0_{\mathrm{cont}}(B, \check{\mathcal E}^\vee(-\Sigma) \otimes L) \, | \, \check{s}(b_j) \in \epsilon_{b_j} \subset \mathcal E^\vee_{b_j} \text{ for all $j \in J$}\}.
\end{multline*}
As this is a closed subspace of $\rH^0_{\mathrm{cont}}(B, \check{\mathcal E}^\vee(-\Sigma) \otimes L)$ this is
 a Banach space.
Moreover,
\begin{equation} 
\label{eqn:semi}
\rH^0_{\mathrm{cont}}(B, \check{\mathcal E}^\vee(-\Sigma) \otimes L)_\Sigma \cap \rH^0_{\mathrm{hol}}(B, \mathcal E^\vee \otimes L) = \rH^0_{\mathrm{hol}}(B, \mathcal E^\vee(-\Sigma) \otimes L).
\end{equation}
The local trivializations of $\mathcal L_d$ over $\Pic^d(B)$ enable us to realize 
\[
\mathcal V^{\mathrm{stop}}_{d, \Sigma} \to \Pic^d(B),
\]
as a locally trivial bundle of Banach spaces.
Let $\tilde{S}^{\mathrm{stop}}_{d, \Sigma} \subset \mathcal V^{\mathrm{stop}}_{d, \Sigma}$ be the open subset parametrizing $(s, L)$ such that $s$ is nowhere vanishing on $B$ 
and $S^{\mathrm{stop}}_{d, \Sigma}$ its quotient by the $\bC^\times$-action \eqref{eqn:c}. 
We have an inclusion
\[
\mathrm{Sect}(\cX/B, d, \Sigma) \hookrightarrow S^{\mathrm{stop}}_{d, \Sigma}.
\]

\subsection*{Comparison}

A class $[s, L] \in S^{\mathrm{stop}}_{d, \Sigma}$ defines a continuous section $\sigma : B \to \mathcal X$ matching $\Sigma$ such that $\sigma^*\mathcal O(1) \approx L$. 

\begin{prop}
\label{prop:comparison}
The continuous map
\[
\Phi_{d, \Sigma} : S^{\mathrm{stop}}_{d, \Sigma} \to \Sect^{\mathrm{stop}}(\mathcal X/B, d, \Sigma) \times \Pic^d(B), \quad  [s, L] \mapsto (\sigma, L)
\]
is a homeomorphism.
\end{prop}

\begin{proof}
Since the line bundle $\mathcal L_d \to B \times \Pic^d(B)$ can be locally topologically trivialized over $\Pic^d(B)$, the map $\Phi_{d, \Sigma}$ is a local homeomorphism over $\Pic^d(B)$. Since $\Phi_{d, \Sigma}$ is bijective, this proves our assertion.
\end{proof}

\section{Comparison of algebraic and semi-topological models}

In this section, we compare (co)homologies of
\[
\mathcal V_{d, \Sigma} \setminus \tilde{S}_{d, \Sigma}, \quad \mathcal V^{\mathrm{stop}}_{d, \Sigma} \setminus \tilde{S}^{\mathrm{stop}}_{d, \Sigma},
\]
following \cite{DT24, DLTT25}.

\subsection*{Stratifications}

We recall several definitions from \cite{DT24}:
\begin{defi}
    A {\em topological poset} is a topological space $P$ and a closed subspace $\mathcal R_P \subset P\times P$ which defines a poset structure on $P$.
 
\end{defi}
\begin{defi}
    Let $P$ be a topological poset and $X$ a topological space. A {\em stratification} of $X$ by $P$ is a closed subset $Z \subset X \times P$ such that for all $p$, the fiber $Z_p$ is a closed subspace of $X$, and $Z_p \supset Z_q$, when $p \leq q$.
\end{defi}
Let 
$$
\mathrm{Hilb}(B)
$$ 
be the Hilbert scheme of $0$-dimensional subschemes on $B$.
We view this space as a  topological poset by inclusion of subschemes.
We define a semi-topological stratification via 
the introduction of 
\[
Z_{d, \Sigma}^{\mathrm{stop}} \subset  \mathcal V^{\mathrm{stop}}_{d, \Sigma} \times \mathrm{Hilb}(B),
\]
the closed subspace whose fiber at $(L, y)$ is the subspace 
\[
(\rH^0_{\mathrm{cont}}(B, \check{\mathcal E}^\vee(-\Sigma) \otimes L)_{\Sigma})_{y}\subset \rH^0_{\mathrm{cont}}(B, \check{\mathcal E}^\vee(-\Sigma) \otimes L)_{\Sigma}
\]
of sections vanishing on the support of $y$:

\ 

\centerline{
\xymatrix@C=0pt{
Z_{d, \Sigma}^{\mathrm{stop}}|_{(L,y)} \subset \ar[d] & Z_{d, \Sigma}^{\mathrm{stop}} \subset &  \mathcal V^{\mathrm{stop}}_{d, \Sigma} \times \mathrm{Hilb}(B) \ar[d] \\
(L,y) \in & & \Pic^d(B)\times \Hilb(B)
}
}

\ 

\noindent
By definition,  
\[
Z_{d, \Sigma}^{\mathrm{stop}}|_{(L, y)} = Z_{d, \Sigma}^{\mathrm{stop}}|_{(L, \mathrm{red}(y))}, 
\]
where $\mathrm{red}(y)$ is the reduced scheme of $y$. 
The inclusion $\mathcal V_{d, \Sigma} \hookrightarrow \mathcal V^{\mathrm{stop}}_{d, \Sigma}$, 
combined with observation \eqref{eqn:semi}, yields 
an algebraic stratification
\[
Z_{d, \Sigma}^{\mathrm{alg}} = \mathcal V_{d, \Sigma} \times_{\mathcal V^{\mathrm{stop}}_{d, \Sigma}} Z_{d, \Sigma}^{\mathrm{stop}}. 
\]

\subsection*{Combinatorial types}
Let $\Sigma = \{\hat{\sigma}\}_{j\in J}$ be a set of admissible $N_j$-th jets, supported in $x_j\in \cX_{b_j}$. 
For an effective $y \in \mathrm{Hilb}(B)$, we express
$$
\deg(y)= \sum_{j\in J} \ell_j +\sum_{i\in I}  m_i,
$$
where $m_i\ge 1$ are multiplicities of $y$ in points $c_i$  outside $\{ b_j\}_{j\in J}$, and $\ell_j$ is the multiplicity of $y$ at $b_j$. In particular, it is possible that $\ell_j =0$.
We define the {\em combinatorial type} of $y$ by setting the multiset 
$$
T(y):=\{ {\boldsymbol{\ell}}; \boldsymbol{m}\}, \quad \quad \boldsymbol{\ell}=\{\ell_j\}_{j\in J}, \quad \boldsymbol{m}=\{m_i\}_{i\in I}.
$$
It is {\em essential}, in the sense of \cite[Definition 3.11]{DT24}, if and only if $y$ is reduced, i.e., all multiplicities are $\le 1$.
In our setting, all combinatorial types are {\em saturated}, in the sense of \cite[Definition 3.3]{DT24}. 
Let 
$$
\mathcal N_T \subset \mathrm{Hilb}(B)
$$ 
be the locally closed subset parametrizing effective divisors of combinatorial type $T$. This defines a stratification into locally closed subsets:
\[
\mathrm{Hilb}(B) = \bigsqcup_{T}\,  \mathcal N_T.
\]
As in \cite[Section 3]{DT24}, 
put 
\begin{equation} 
\label{eqn:gamma}
\gamma(y):= (n + 1)\left(\sum_{i \in I} m_i \right) + \sum_{j\in J}(\ell_j + n\max\{\ell_j - N_j - 1, 0\}),
\end{equation} 
When $y$ is reduced, 
this is the expected codimension of the incidence condition imposed by $y$, i.e., the expected codimension of 
$$
Z^{\mathrm{alg}}_{d, \Sigma}|_{(L, y)}\subset \mathcal V_{d, \Sigma}|_L.
$$
We note that $\mathrm{rank}(y)$ from \cite[Example 3.15]{DT24} equals $\deg(y)$, in our situation.

\subsection*{Semi-algebraic approximation}

Following \cite[Section 6.3]{DT24}, let $M$ be a very ample line bundle on $B$ and $\overline{M}$ its antiholomorphic bundle. Let 
$$
\cW_k \subset \mathcal V^{\mathrm{stop}}_{d, \Sigma}
$$
be the subbundle such that its fiber $(\cW_k)_L$ over $L$ is the image of
\[
\rH^0(B, \mathcal E^\vee(-\Sigma) \otimes L \otimes M^k) \otimes \rH^0_{\mathrm{anti}}(B, \overline{M}^k) \to \rH^0_{\mathrm{cont}}(B, \check{\mathcal E}^\vee(-\Sigma)\otimes L)_{\Sigma},
\]
where the map is given by a natural topological trivialization 
$$
\mathbb C \times B \cong M \otimes \overline{M}.
$$
This map is injective, by \cite[Lemma 6.17]{DT24}, thus its image has dimension independent of $L$. Hence $\cW_k$ is a finite-dimensional semi-algebraic vector bundle over $\Pic^d(B)$, and we have inclusions
\[
\mathcal V_{d, \Sigma} = \cW_0 \subset \cW_1 \subset \cdots \subset \cW_k \subset \cdots \subset \mathcal V^{\mathrm{stop}}_{d, \Sigma}, 
\quad k\in \bZ_{\ge 0}.
\]
Indeed, the inclusions follow from the fact that $M$ is very ample and we have the multiplication map
\begin{align*}
&\rH^0(B, \mathcal E^\vee(-\Sigma) \otimes L \otimes M^k) \otimes \rH^0_{\mathrm{anti}}(B, \overline{M}^k)\\ & \qquad \qquad \qquad\to
\rH^0(B, \mathcal E^\vee(-\Sigma) \otimes L \otimes M^{k+1}) \otimes \rH^0_{\mathrm{anti}}(B, \overline{M}^{k+1}) 
\end{align*}
by a section $s \otimes \bar{s}$, where $s \in \rH^0(B, M)$.
Using the argument of \cite[Lemma 7.2]{Aumonier}, one can prove that $\cup_k \cW_k$ is dense in $\mathcal V^{\mathrm{stop}}_{d, \Sigma}$.
Moreover,

\centerline{
\xymatrix@C=0pt{
Z_{d,\Sigma,k} :=\cW_k \times_{\mathcal V^{\mathrm{stop}}_{d, \Sigma}} Z_{d, \Sigma}^{\mathrm{stop}} \subset &  \cW_k \times \mathrm{Hilb}(B) \ar[d] \\
& \Pic^d(B)\times \Hilb(B)
}
}

\noindent
is also a semi-algebraic stratification of $\cW_k$, 
as the condition on $(\cW_k)|_L$ imposed by $y \in \mathrm{Hilb}(B)$ is linear.

\subsection*{The bar complex}
Put 
$$
Z_k:=Z_{d,\Sigma,k},
$$ 
and let $U \subset \Pic^d(B)$ be an open subset, in the Euclidean topology. Viewing $Z_k$ as a bundle over $\Pic^d(B)$, we 
let $Z_{k, U}$ be its restriction to $U$. 
Consider
$$
\rP:=\cup_T \cN_T \subset \mathrm{Hilb}(B),
$$ 
a downward closed proper union, over finitely many $T$. 
Following \cite[Section 5]{DT24}, the {\em bar complex}
 $$
 \rB(\rP, Z_{k, U})
 $$ 
is a simplicial space with $r$-simplices
\[
\{ (L, y_0 < \cdots < y_r, s) \, |\, L \in U, \quad  y_i \in \rP, \quad s \in (Z_{k, U})_{(L, y_r)} )\}.
\]
We denote its geometric realization by $\mathbf{B}(\rP, Z_{k, U})$.

\subsection*{Unobstructedness}

\begin{prop}
\label{prop:unobstructedness}
Let $y \in \mathrm{Hilb}(B)$ be reduced, of combinatorial type 
$$
T(y):=\{ {\boldsymbol{\ell}}; \boldsymbol{m}\}, \quad \quad \boldsymbol{\ell}=\{\ell_j\}_{j\in J}, \quad \boldsymbol{m}=\{m_i\}_{i\in I}.
$$
There exists a constant $A({\mathcal E, \deg(\Sigma)})$ 
such
that for 
$$
|I| \leq d -A({\mathcal E, \deg(\Sigma)})
$$
the real codimension of 
$$
(Z_{d,\Sigma,k})_{(L, y)}=(Z_{k})_{(L, y)} \subset (\cW_{k})_L 
$$
is equal to $2\gamma(y)$, defined in \eqref{eqn:gamma}. 
\end{prop}

\begin{proof}
Since $y$ imposes a linear condition, the real codimension 
of 
\[
Z_{k}|_{(L, y)} \subset (\cW_k)_L
\]
is less than or equal $2\gamma(y)$. 
Equality holds if the real codimension of 
$$
(Z_{d, \Sigma}^{\mathrm{alg}})_{(L, y)} \subset (\mathcal V_{d, \Sigma})_L 
$$
is equal to $2\gamma(y)$, which we now prove. 
For $d\gg 0$, depending on $\mathcal E, \deg(\Sigma$), we have a surjection
\[
\rH^0(B, \mathcal E^\vee\otimes L) \to \oplus_{j \in J} \mathcal E^\vee|_{y_j}, \quad y_j = (N_j + 1)b_j, \quad \forall j\in J, 
\]
and the conditions imposed by $b_j \in \mathrm{Supp}(y)$ are independent. 
There is a constant  $A(\cE,\deg(\Sigma))$ such that for 
$$
d - |I|\geq A(\cE,\deg(\Sigma)),
$$
the homomorphism
\[
\rH^0(B, \mathcal E^\vee\otimes L )\to \oplus_{j \in J} \mathcal E^\vee|_{y_j} \oplus \oplus_{i \in I} \mathcal E^\vee_{c_i},
\]
is surjective. 
In that range, 
$$
(Z_{d, \Sigma}^{\mathrm{alg}})_{(L, y)} \subset (\mathcal V_{d, \Sigma})_L 
$$
has the expected real codimension $2\gamma(y)$.
\end{proof}

\subsection*{Approximation}

We show that 
\begin{equation}
\label{eqution:kappavsdegree}
   \deg (T)  - 2|\Sigma| \leq \kappa(T):= 2\gamma(T) - \deg (T) - 2|\mathrm{Supp}(T)|  
\end{equation}
Indeed, using (\ref{eqn:gamma}) we have
\begin{align*}
  \kappa(T) &= (2n+1)\sum_{i \in I} m_i  -2|I|\\
  & \qquad\qquad + \sum_{j\in J}(\ell_j + 2n\max\{\ell_j - N_j - 1, 0\} -2\min\{\ell_j, 1\})\\
  &\geq \sum_{i \in I} m_i + \sum_{j \in J} \ell_j -2|\Sigma| = \deg(T) - 2|\Sigma|.
\end{align*}
Fix a positive integer $R$ and define
\[
\rP = \{ y \in \mathrm{Hilb}(B) \, |\, \deg (y) \leq R + 2|\Sigma| \},
\]
a downward closed proper union of finitely many $\mathcal N_T$ with $\deg(T) \leq R + 2|\Sigma|$.
One of the main applications of \cite[Theorem 5.9]{DT24} is:

\begin{prop}
\label{prop:barcomplex}
Suppose that 
$$
d \geq R  + 2|\Sigma|+ A({\mathcal E, \deg(\Sigma)})+1.
$$
Then the map 
$$
\rB(\rP, Z_{k, U}) \to \mathrm{im}(Z_{k, U}|_{\rP} \to \mathcal W_k|_U)
$$ induces a homomorphism
\[
\rH^i_c(\mathrm{im}(Z_{k, U}|_{\rP} \to \mathcal W_k|_U), \mathbb Z)\to \rH^i_c(\mathbf{B}(\rP, Z_{k, U}), \mathbb Z),
\]
which is an isomorphism when $i > \dim(\cW_k) - R - 2$ and a surjection when $i = \dim(\cW_k) - R - 2$, where $\dim(\cW_k)$ is the real dimension of the semi-algebraic bundle $\cW_k$.
\end{prop}

\begin{proof}
This follows from a version of \cite[Theorem 5.9]{DT24}. To verify the assumptions, observe first that
$\rP$ is downward closed and proper. Furthermore, 
for $y \in \rP$ and $y \prec y'$ such that $y'$ is reduced, 
we have 
$$
|I(y')|\leq \deg(y') \leq R  + 2|\Sigma| + 1.
$$
By Proposition~\ref{prop:unobstructedness}, 
$$
(Z_{k})_{L, y} \subset (\cW_{k})_L
$$
has expected real codimension $2\gamma(y)$.
Finally, by (\ref{eqution:kappavsdegree}), $\kappa(T)\leq R$ implies that $T$ is the type of $\rP$. 
\end{proof}

\subsection*{The main result}


\begin{theo}
\label{theo:projectivebundlealgvstop}
    Assume that 
    $$
    d \geq R  + 2|\Sigma| +  A({\mathcal E, \deg(\Sigma)}) + 1.
    $$
    Then the inclusion
    \[
    \Sect(\mathcal X/B, d, \Sigma) \hookrightarrow S^{\mathrm{stop}}_{d, \Sigma},
    \]
    is homology $R$-connected, i.e., 
   the induced homomorphism
    \[
    \rH_i(\Sect(\mathcal X/B, d, \Sigma) , \mathbb Z) \to \rH_i(S^{\mathrm{stop}}_{d, \Sigma}, \mathbb Z),
    \]
    is an isomorphism when $i < R$ and an injection when $i = R$.
\end{theo}

\begin{proof}
    We follow the proof of \cite[Theorem 7.1]{DT24}:

\

\noindent
{\em Step 1.}
    Both $\tilde{S}_{d, \Sigma} \to \Sect(\mathcal X/B, d, \Sigma)$ and $\tilde{S}^{\mathrm{stop}}_{d, \Sigma}\to S^{\mathrm{stop}}_{d, \Sigma}$ are $\mathbb C^\times$-torsors,
    so they induce the following diagram:
    \begin{equation}
    \xymatrix{
    \Sect(\mathcal X/B, d, \Sigma) \ar@{^{(}->}[r]\ar[d] & S^{\mathrm{stop}}_{d, \Sigma} \ar[d]\\
   B\mathbb C^\times \ar@{=}[r] & B \mathbb C^\times.
    }
    \end{equation}
    Applying the Leray spectral sequence to both vertical maps, the theorem follows from homology $R$-connectedness of   
    the inclusion
    \[
    \tilde{S}_{d, \Sigma} \hookrightarrow \tilde{S}^{\mathrm{stop}}_{d, \Sigma}.
    \]

\

\noindent
{\em Step 2.}
Since the inclusion is a continuous map, compatible with projections to $\Pic^d(B)$, 
it follows from the two Leray spectral sequences over $\Pic^d(B)$ that it suffices to show that for a basis of open subsets $U \subset \Pic^d(B)$, the restriction
    \[
    \tilde{S}_{d, \Sigma}|_U \hookrightarrow \tilde{S}^{\mathrm{stop}}_{d, \Sigma}|_U,
    \]
    is homology $R$-connected.

\

\noindent
{\em Step 3.}
We have realized 
$$
\tilde{S}^{\mathrm{stop}}_{d, \Sigma}|_U\subset \mathcal V^{\mathrm{stop}}_{d, \Sigma}|_U
$$
as an open subset of a Banach bundle over $U$. Taking a sufficiently small $U$, we trivialize this bundle. By \cite[Proposition 6.16]{DT24}, the inclusion 
    \[
    \cup_k(\tilde{S}^{\mathrm{stop}}_{d, \Sigma}|_U \cap (\cW_k|_U)) \hookrightarrow \tilde{S}^{\mathrm{stop}}_{d, \Sigma}|_U
    \]
    is a weak homotopy equivalence. 
    Thus, it suffices to show that
    \[
    \tilde{S}_{d, \Sigma}|_U \hookrightarrow \tilde{S}_{k}|_U :=\tilde{S}^{\mathrm{stop}}_{d, \Sigma}|_U \cap (\cW_k|_U),
    \]
    is homology $R$-connected, for sufficiently large $k$.

    \

\noindent
{\em Step 4.}
Since $\tilde{S}_{d, \Sigma}|_U = \tilde{S}_{k}|_U \cap \mathcal V_{d, \Sigma}|_U$, by \cite[Proposition 2.2]{DT24}, the induced homomorphism
\[
\rH_i(\tilde{S}_{d, \Sigma}|_U, \mathbb Z) \to \rH_i(\tilde{S}_{k}|_U , \mathbb Z),
\]
is Poincar\'e dual to the Gysin map
\[
\rH_c^{2\dim(\mathcal V_{d, \Sigma}) -i}(\tilde{S}_{d, \Sigma}|_U, \mathbb Z) \to \rH_c^{\dim(\cW_k) -i}(\tilde{S}_{k}|_U , \mathbb Z).
\]
It suffices to prove that this Gysin map is an isomorphism, when $i < R$, and a surjection, when $i = R$.

\

\noindent
{\em Step 5.}
Let
\[
C^{\mathrm{alg}} = \mathcal V_{d, \Sigma}|_U \setminus \tilde{S}_{d, \Sigma}|_U, \quad C_k = \cW_k|_U \setminus \tilde{S}_{k}|_U.
\]
We have a commuting diagram of long exact sequences of cohomology with compact supports:

\ 

\centerline{
    \xymatrix@C=0.6em{
 & \rH_c^{2\dim(\mathcal V_{d, \Sigma}) -i}(\tilde{S}_{d, \Sigma}|_U, \mathbb Z) \ar[r]\ar[d] & \rH_c^{2\dim(\mathcal V_{d, \Sigma}) -i}(\mathcal V_{d, \Sigma}|_U, \mathbb Z) \ar[r]\ar[d] & \rH_c^{2\dim(\mathcal V_{d, \Sigma}) -i}(C^{\mathrm{alg}}, \mathbb Z)  \ar[d]^{\psi_i}&  \\
 & \rH_c^{\dim(\mathcal W_{k}) -i}(\tilde{S}_{k}|_U, \mathbb Z) \ar[r] & \rH_c^{\dim(\mathcal W_{k}) -i}(\mathcal W_{k}|_U, \mathbb Z) \ar[r] & \rH_c^{\dim(\mathcal W_{k}) -i}(C_k, \mathbb Z).  &  
    }
}

\

Here the vertical maps are the Gysin maps.
It suffices to show that $\psi_i$ 
is an isomorphism when $i < R+1$ and a surjection when $i = R+1$.

\

\noindent
{\em Step 6.}
By Proposition~\ref{prop:barcomplex}, it suffices to show that the Gysin map
\[
\rH_c^{2\dim(\mathcal V_{d, \Sigma}) -i}(\mathbf B(\rP, Z^{\mathrm{alg}}|_U), \mathbb Z) \to
\rH_c^{\dim(\mathcal W_{k}) -i}(\mathbf B(\rP, Z_{k, U}), \mathbb Z)
\]
is an isomorphism when $i < R+1$ and a surjection when $i = R+1$.
In turn, this follows from a version of \cite[Theorem 5.6]{DT24}, in our setting.

\end{proof}

\ 

\section{Homological stability}
\label{sec:homological}

\subsection*{Projective bundles over curves}

\begin{theo}
\label{theo:homologicalstabilityI}
Let
\[
\pi : \mathcal X = \mathbb P(\mathcal E) \to B
\]
be the projectivization of a vector bundle, of relative dimension $n\ge 1$.
Let $\Sigma,\Sigma'$ be non-empty sets of admissible jets for $\pi$, such that 
$$
\pi(\Sigma) = \pi(\Sigma')\quad \text{ and } \quad \deg(\Sigma)\ge \deg(\Sigma').
$$
Let 
\[
\ell (d) = d - 2|\Sigma|-A({\mathcal E, \deg(\Sigma)}) -2,
\]
where $A(\mathcal E, \deg(\Sigma))$ is the constant from Proposition~\ref{prop:unobstructedness}.
Then, for all $i \leq \ell(d)$, we have isomorphisms
\begin{align*} 
\rH_i(\Sect(\mathcal X/B, d, \Sigma), \mathbb Z) &  \cong \rH_i(\Sect(\mathcal X/B, d + 1, \Sigma), \mathbb Z) \\
& \cong \rH_i(\Sect(\mathcal X/B, d, \Sigma'), \mathbb Z).
\end{align*}
\end{theo}

\begin{proof}
The first isomorphism follows from Theorem~\ref{theo:projectivebundlealgvstop}, Proposition~\ref{prop:comparison}, and Corollary~\ref{coro:homologicalstabilitytop}, with $\beta$ being the class of a vertical line. Indeed, we have
\begin{align*}
    \rH_i(\Sect(\mathcal X/B, d, \Sigma), \mathbb Z) & \cong \rH_i(S^{\mathrm{stop}}_{d, \Sigma}, \mathbb Z)\\
    & \cong \rH_i(\Sect^{\mathrm{stop}}(\mathcal X/B, d, \Sigma) \times \Pic^d(B), \mathbb Z) \\
    & \cong \rH_i(\Sect^{\mathrm{stop}}(\mathcal X/B, d + 1, \Sigma) \times \Pic^d(B), \mathbb Z)\\
    & \cong \rH_i(S^{\mathrm{stop}}_{d + 1, \Sigma}, \mathbb Z)\\
    & \cong \rH_i(\Sect(\mathcal X/B, d + 1, \Sigma), \mathbb Z).
\end{align*}
The second isomorphism follows from Theorem~\ref{theo:projectivebundlealgvstop}, Proposition~\ref{prop:comparison}, Proposition~\ref{prop:pointsindependnetgeneral}, and Proposition~\ref{prop:independentofjets}. We have
\begin{align*}
    \rH_i(\Sect(\mathcal X/B, d, \Sigma), \mathbb Z) & \cong \rH_i(S^{\mathrm{stop}}_{d, \Sigma}, \mathbb Z)\\
    & \cong \rH_i(\Sect^{\mathrm{stop}}(\mathcal X/B, d, \Sigma) \times \Pic^d(B), \mathbb Z) \\
    & \cong \rH_i(\Sect^{\mathrm{stop}}(\mathcal X/B, d , \Sigma') \times \Pic^d(B), \mathbb Z)\\
    & \cong \rH_i(S^{\mathrm{stop}}_{d, \Sigma'}, \mathbb Z)\\
    & \cong \rH_i(\Sect(\mathcal X/B, d, \Sigma'), \mathbb Z).
\end{align*}
\end{proof}

\subsection*{Conic bundles over curves}
Let
\[
\pi : \mathcal S \to B,
\]
be a smooth conic bundle over $B$, i.e., $\mathcal S$ is a smooth projective surface such that $\omega_{\pi}^{-1}$ 
is $\pi$-relatively ample.
Any singular fiber is the union of two lines meeting at a point.
Let $\mathfrak d$ be the number of singular fibers.
There are $2^{\mathfrak d}$ birational morphisms 

\

\centerline{
\xymatrix{ 
\mathcal S\ar[r]^(.3){\varphi_r}\ar[d]_{\pi} & \mathcal S_r= \mathbb P(\mathcal E_r) \ar[d]^{\pi_r}  \\ 
    B    \ar@{=}[r] &         B
}
}

 \
 
\noindent 
contracting vertical $(-1)$-curves;  here $\mathcal E_r$ is normalized so that
\[
\omega^{-1}_{\pi_r}\cong \mathcal O_{\mathbb P(\mathcal E_r)/B}(2)\otimes \mathcal O(\mathcal S_b)^{\otimes \varepsilon(r)},
\]
where $\mathcal S_b$ is a general fiber of $\pi$ and $\varepsilon(r) = 0$, or $1$. 
Given a section $\sigma$ of $\pi$ class $\alpha$ we let $\varphi_{r(\alpha)}$ be
the morphism contracting curves meeting $\sigma$.  
Put 
\begin{equation} 
\label{eqn:deg}
h=\deg(\sigma^*\omega^{-1}_\pi),
\quad 
d = \deg(\sigma^*\varphi_{r(\alpha)}^{*}\mathcal O_{\mathbb P(\mathcal E_{r(\alpha)})/B}(1))
\end{equation}
In particular,
\[
h = 2d + \epsilon(r) -\mathfrak d.
\]

Let $\Sigma$ be a set of admissible jets for $\pi$ and $\Sigma_{r(\alpha)}$ 
the set of admissible jets for $\pi_{r(\alpha)}$ such that $\Sigma$ is the strict transform of $\Sigma_{r(\alpha)}$. 
We have
\[
\deg(\Sigma_{r(\alpha)}) = \deg(\Sigma) + \mathfrak d, \quad |\Sigma_{r(\alpha)}| \leq |\Sigma| + \mathfrak d.
\]

Proposition~\ref{prop:blowupmodel} implies that 
there is an isomorphism
\begin{equation}
     \label{eqn:birational contraction}
\Sect(\mathcal S/B,h,\Sigma)\simeq 
\Sect(\mathcal S_{r(\alpha)}/B,h + \mathfrak d ,\Sigma_{r(\alpha)}).
\end{equation}

\begin{theo}
Let 
$$
\pi: \mathcal S\to B
$$ 
be a smooth conic bundle. Let $\Sigma,\Sigma'$ be non-empty sets of admissible jets for $\pi$, such that 
$$
\pi(\Sigma) = \pi(\Sigma')\quad \text{ and } \quad \deg(\Sigma)\ge \deg(\Sigma').
$$
Let $\alpha$ be the class of a section of $\pi$, $h$ its degree as in \eqref{eqn:deg},  and 
\[
\ell (h) := \frac{h}{2}  -\frac{1}{2}- 2|\Sigma|-\frac{\mathfrak d}{2}-A(\mathcal E_{r(\alpha)}, \deg(\Sigma)+\mathfrak d) -2.
\]
Then, for all $i \leq \ell(h)$, we have isomorphisms
\begin{align*} 
\rH_i(\Sect(\mathcal S/B, h, \Sigma), \mathbb Z) &  \cong \rH_i(\Sect(\mathcal S/B, h + 2, \Sigma), \mathbb Z) \\ 
& \cong \rH_i(\Sect(\mathcal S/B, h, \Sigma'), \mathbb Z).
\end{align*}
\end{theo}

\begin{proof}
    This follows from \eqref{eqn:birational contraction} and Theorem~\ref{theo:homologicalstabilityI}.
\end{proof}

\subsection*{Quadric surface bundles}
\label{sect:gen}

Let 
$$
\pi:\cX\to B
$$ 
be a smooth quadric surface bundle with relative Picard rank one and singular fibers with at most one singularity, of type $\mathsf A_1$. 
Let 
$$
F_1(\cX)\to D\stackrel{\iota}{\lra} B
$$
be the Stein factorization of the map
from the space of lines in the fibers of $\pi$ to $B$ (see, e.g., \cite[Section 3]{HT12} for more details).  
The covering involution $\iota$ is branched along the discriminant divisor $\mathfrak d$ of $\pi$; the map 
$F_1(\cX)\to D$ is a {\em smooth} $\bP^1$-bundle, i.e., 
$$
\pi_{\cY}: \cY:=F_1(\cX)=\bP(\cE)\to D,
$$
for a rank-2 vector bundle $\cE$ on $D$. Every point $x\in\cX$, in a smooth fiber of $\pi$, gives rise to points $y,y'\in \cY$, in distinct fibers of $\pi_{\cY}$. 
In detail, 
let
\[
\mathcal D \subset \mathcal X\times_B F_1(\mathcal X).
\]
be the universal family of lines.
Given a section $\sigma : B \to \mathcal X$, 
define $\mathcal D_\sigma$ as the fiber product
\begin{equation*}
    \xymatrix{
\mathcal D_\sigma \ar@{^{(}->}[r] \ar[d]& B \times_B F_1(\mathcal X) \ar[d]^{\sigma \times \mathrm{id}} \\
\mathcal D \ar@{^{(}->}[r] & \mathcal X\times_B F_1(\mathcal X)  .  
    }
    \end{equation*}
    Then 
    \[
    \mathcal D_\sigma \to B \times_B F_1(\mathcal X) \to B
    \]
    is a double cover ramified along $\mathfrak d$, so that $\mathcal D_\sigma$ is isomorphic to $D$. We obtain a section 
    \[
    \tau : D \to \mathcal F_1(\mathcal X).
    \]
    This construction applies to jets as well:
    let 
    \[
    \hat{\sigma} : \Spec (\mathbb C[z]/(z^{N+1})) \to \mathcal X
    \]
    be an admissible $N$th-jet. This yields an admissible jet
    \[
    \hat{\tau}: \mathcal D_{\hat{\sigma}} \to F_1(\mathcal X).
    \]
    When $\hat{\sigma}$ is supported in a smooth fiber, 
    \[
    \mathcal D_{\hat{\sigma}} \to \Spec (\mathbb C[z]/(z^{N+1}))
    \]
    is \'etale, so it consists of two copies of $\Spec (\mathbb C[z]/(z^{N+1}))$.
    When $\hat{\sigma}$ is supported in a singular fiber, $\mathcal D_{\hat{\sigma}}$ is isomorphic to $\Spec (\mathbb C[z]/(z^{2N+2}))$.
    In this way, a set of admissible jets $\Sigma_B$ induces a set of admissible jets $\Sigma_D$, with 
    \[
    \deg(\Sigma_D) = 2\deg(\Sigma_B), \quad |\Sigma_D| \leq 2|\Sigma_B|.
    \]
\begin{lemm} 
\label{lemm:bi}
Let $\Sigma_B$ be a set of admissible jets for $\pi : \mathcal X \to B$
and $\Sigma_D$ the induced set of admissible jets of $\pi_{\cY} : \mathcal Y \to D$.
There is an isomorphism
$$
\Sect(\cX/B,h,\Sigma_B)\simeq 
\Sect(\cY/D, d,\Sigma_D),
$$
for
$$
h = 2d + \deg(\mathcal E) - \frac{|\mathfrak d|}{2}.
$$
\end{lemm}

\begin{proof}
    The last formula follows from the normalization of $\mathcal E$ in \cite[Section 3]{HT12}.
\end{proof}

\begin{theo}
\label{theo:homologicalstabilityIII}
Let 
$$
\pi:\cX\to B
$$ 
be a smooth quadric surface bundle with relative Picard rank one and singular fibers with at most one $\mathsf A_1$-singularity. 
Let $\Sigma_B, \Sigma_B'$ be non-empty sets of admissible jets for $\pi$, such that 
$$
\pi(\Sigma_B) = \pi(\Sigma_B')\quad \text{ and } \quad \deg(\Sigma_B)\ge \deg(\Sigma_B').
$$
Let 
\[
\ell (h) = \frac{h}{2} -\frac{\deg(\mathcal E)}{2} + \frac{|\mathfrak d|}{4} - 4|\Sigma_B|-A(\mathcal E, 2\deg(\Sigma_B)) -2.
\]
Then, for all $i \leq \ell(h)$, we have isomorphisms
\begin{align*}
\rH_i(\Sect(\mathcal X/B, h, \Sigma_B), \mathbb Z)  & \cong \rH_i(\Sect(\mathcal X/B, h + 2, \Sigma_B), \mathbb Z)\\ 
&\cong  \rH_i(\Sect(\mathcal X/B, h, \Sigma_B'), \mathbb Z).
\end{align*}
\end{theo}

\begin{proof}
    This follows from Lemma~\ref{lemm:bi} and Theorem~\ref{theo:homologicalstabilityI}.
\end{proof}

\bibliographystyle{alpha}
\bibliography{quadricsurfacebundle}

\end{document}